\documentclass[11pt,a4paper,notitlepage,final,onecolumn,oneside]{article}
\usepackage{exscale}
\usepackage{latexsym}
\usepackage{calc}
\usepackage{amssymb}

\begin{document}

\newcommand{\bks}{\backslash}


\title{The polynomial algorithm for graphs' isomorphism testing}
\author{Aleksandr Golubchik}

\maketitle

\begin{abstract}
A polynomial algorithm for graphs' isomorphism testing will be
constructed in assumption that there exists a corresponding polynomial
algorithm for graphs with trivial automorphism group.
\end{abstract}


In this paper we use the term ``graph'' for any partition of the cartesian
square $V^2$ of the set of $n$ vertices $V$. In other words, the graph
$G(V)$ denotes a color digraph with color vertices.

It is known \cite{Weis} that the graphs' isomorphism problem is equivalent
to the problem of determination of orbits of graph's automorphism group.
We will construct the polynomial algorithm for latter problem.


\section{Lemma}
Let $A, B$ be permutation groups on $V$, and $Orb(A), Orb(B)$ be
corresponding systems of orbits on $V$ then a partition $Q=Orb(A)\cup
Orb(B)$ is the system of orbits $Orb(C)$ of group $C$ generated by
subgroups $A, B$ (where the union of partitions means the union of its
classes containing non-empty intersection).

\paragraph*{\em Proof}
\begin{enumerate}
\item
As any permutation $a\in A$ keeps partition $Orb(A)$ and any permutation
$b\in B$ keeps partition $Orb(B)$ then any product of permutations from
$A, B$ keeps partition $Q$.
\item
Let $U$ be a class of $Q$ and points $x,y\in U$ then there exists a
sequence of classes $V_1,W_1,V_2,W_2,V_3,\ldots$ so that $V_i\in Orb(A)$,
$W_i\in Orb(B)$, any two of nearby classes have non-empty intersection and
points $x,y$ are located in the utmost classes. It follows that there
exists a product of elements from $A,B$ that transforms $x$ in $y$.
\end{enumerate}


\section{Algorithm}
Let $A$ be a polynomial algorithm that colors vertices of the graph $G(V)$
in ordered colors and thus determines an ordered partition on the set $V$.
Let a graph $AG$ obtained from graph $G$ by action of algorithm $A$ has
the same automorphism group as graph $G$. Let also this algorithm
recognize graphs with trivial automorphism group so that it colors every
vertex of such graph in its unique color.

Now we can construct a polynomial algorithm that finds the orbits
of automorphism group of any graph $G(V)$.

\begin{enumerate}
\item
Obtain graphs $R=AG$, $R(x_1)=AG(x_1)$, $R(x_1,x_2)=AG(x_1,x_2)$,
\ldots, $R(x_1,x_2,\ldots,x_k)=AG(x_1,x_2,\ldots,x_k)$ consecutively, by a
consecutive coloring of vertices $x_1,\ldots,x_k$ from the set $V$, until
the graph $R(x_1,\ldots,x_k)$ with regular automorphism group from the
initial graph $G$ is obtained.

It means that the graph $R(x_1,\ldots,x_k)$ has a non-trivial automorphism
group, but for any other vertex $x_{k+1}$ a graph $G(x_1,\ldots,x_{k+1})$
obtained from the graph $R(x_1,\ldots,x_k)$ by coloring a vertex $x_{k+1}$
has the trivial automorphism group.

Of course, if $G$ has trivial automorphism group, the algorithm stops
immediately.

\item
By testing the graphs $R(x_1,\ldots,x_{k+1})=AG(x_1,\ldots,x_{k+1})$ on
isomorphism, construct the system of orbits of the group
$Aut(G(x_1,\ldots,x_k))$.
\end{enumerate}

It is possible to iterate steps 1 and 2 taking every time a new sequence
of fix-vertices so as to obtain a new {\it automorphic partition} (a
system of orbits of some subgroup of $Aut(G)$) on $V$.  By union of
automorphic partitions we will every time (accordingly to Lemma) obtain a
larger and larger automorphic partition until we will have obtained a
system of orbits of automorphism group of initial graph $G$.

We can notice that the larger automorphic partition is the more restricted
is the possibility for the choice of fix-vertices, because all vertices in
the same orbit are identical for such choice.

By finding the isomorphic graphs we find not only automorphic partitions,
but also automorphisms. So by constructing the system of orbits of
$Aut(G)$ we also obtain a defining system of this group.


\section{How to select fix-vertices}
At first we can see that it is sufficient to find at most $n-1$
automorphic partitions, having non-trivial intersection, in order to
obtain the system of orbits of $Aut(G)$ as union of these partitions. More
precise, if a system of orbits $Q$ of some subgroup of $Aut(G)$ is already
found then it is sufficient to discover further at most $|Q|-1$
corresponding automorphic partitions.

It can be used different strategies for separation of fix-vertices.
One of these can be based on a partial order of set $V$ that orders
vertices by its number of previous fixations.

The main problem is to verify whether an obtained system of orbits $Q$ is
$Orb(Aut(G))$. Let $O_1,O_2\in Q$. Now we will give the verifying algorithm
that verifies whether there exists an automorphism of graph $G$ connecting
these two orbits (suborbits). It is sufficient to search a corresponding
automorphism for any pair of vertices $o_1\in O_1$ and $o_2\in O_2$.

\begin{enumerate}
\item
By fixation of vertices, obtain the graph
$R(x_1,\ldots,x_k)=AG(x_1,\ldots,x_k)$ relative to that the vertices
$o_1,o_2$ are equivalent (have the same color), but for any other vertex
$x_{k+1}$ the graph $R(x_1,\ldots,x_{k+1})=AG(x_1,\ldots,x_{k+1})$
distinguishes the vertices $o_1$ and $o_2$.

\item
Obtain the graphs $T(o_1)\equiv R(x_1,\ldots,x_k,o_1)$ and $T(o_2)\equiv
R(x_1,\ldots,x_k,o_2)$.

If connecting automorphism for $o_1$ and $o_2$ exists then
the graphs $T(o_1)$ and $T(o_2)$ are isomorphic and, accordingly to the
result of the previous step, have the empty intersection of the same color
classes.

\item
Using the whole algorithm from beginning, find a system of orbits of
automorphism groups of graphs $T(o_1)$ and $T(o_2)$ and examine these
graphs on isomorphism.

\end{enumerate}

It is clear that, in order to finish the verifying algorithm, we need to
examine some tree of graphs.

It is clear that the number of levels of consecutive pairs of vertices
(the same as $o_1,o_2$) of such tree  is smaller than $\log_2n$, because
the graphs $T(o_1)$, $T(o_2)$ halve the set $V\bks \{x_1,\ldots, x_k\}$.
And hence the total number of graphs of this tree is smaller than $n$.

This consideration is sufficient for polynomiality of whole algorithm.

There is one more property of the base tree of whole algorithm that
ascertains the polynomiality of the described algorithm: the total number
of graphs of the base tree, that need to be investigated, has a
non-principal difference for various cases and is by order not greater
than $n^2$, because

\begin{itemize}
\item
the greater the number of levels of the tree, the smaller the number of
branches of every level (as for complete graph in the limit, that requires
algorithm $A$ for about $n^2$ graphs )

\end{itemize}
and vice versa,
\begin{itemize}
\item
the greater the number of branches that levels of tree have, the
smaller the number of levels of the tree (as for graph with regular
automorphism group in the limit, that gives $n$ graphs for testing of
their pairs on isomorphism).
\end{itemize}

The complexity of algorithm $A$ is obviously greater than the complexity
of graphs isomorphism testing for graphs with trivial automorphism group.
Hence the case ``more levels'' has greater complexity than the case ``more
branches''. It follows that the whole algorithm has its greatest complexity
for complete graph and this complexity is equal $n^2|A|$, where $|A|$ is
complexity of algorithm $A$.


\section{Conclusion}


\subsection{Algorithm $A$}
The sequence of graph stabilization algorithms can be constructed in the
way B.~Weisfeiler described in \cite{Weis}. We will denote such algorithms
as $A_k$. The algorithm $A_k$ stabilizes the graph $G$ by stabilizing the
structure on $V^k$, generated by this graph. B.~Weisfeiler examined the
algorithm $A_2$ and began to examine the algorithm $A_3$.

The given algorithm puts up a question whether there exists a natural $l$
so that for all $k\geq l$ algorithm $A_k$ recognizes the graphs with
trivial automorphism group. To author an existence of counter-example for
$l=3$ is unknown, i.e. an existence of graph with trivial automorphism
group in that any two isomorphic triangles have the same number of
spanning isomorphic quadrangles.


\subsection{Perspectives}
The algorithms $A_k$ present a combinatorial direction of investigation of
graphs' symmetries. This direction for $k>2$ cannot have a two-place
algebraic interpretation. It is obviously more difficult but substantially
stronger for considered problem than a generalization of $A_2$ on the
algorithms acting on $V^2$ that leads to conception of distance-regular
graph systematically described in \cite{Bannai}.

On the way to solve the graphs isomorphism problem Author has discovered
the original combinatorial objects, that was, most likely, not investigated
earlier. These combinatorial objects present the transformation of
principal properties of graph symmetries on the set $C_k$ of $k+1$
classes from $V^k$ that leads to conception of being $C_k$ assembled or
non-assembled in the subset of $V^{k+1}$.  The simplest example of this
representation on $V^2$ we obtain in the next matrix notation:
\begin{enumerate}
\item Assembled case:
\begin{displaymath}
C'_2=
\left\{
\left[
\begin{array}{c}
  12\\
  45
\end{array}
\right],
\left[
\begin{array}{c}
  23\\
  56
\end{array}
\right],
\left[
\begin{array}{c}
  31\\
  64
\end{array}
\right]
\right\}
=
\left[
\begin{array}{c}
  123\\
  456
\end{array}
\right]
\end{displaymath}

\item Non-assembled case:
\begin{displaymath}
C''_2=
\left\{
\left[
\begin{array}{c}
  12\\
  45
\end{array}
\right],
\left[
\begin{array}{c}
  23\\
  56
\end{array}
\right],
\left[
\begin{array}{c}
  34\\
  61
\end{array}
\right]
\right\}
\neq
\left[
\begin{array}{c}
  123\\
  456
\end{array}
\right]
\end{displaymath}
\end{enumerate}

Both sets $C'$ and $C''$ have the same projection on the set $V$ that
consists of subsets
\begin{displaymath}
\left[
\begin{array}{c}
  1\\
  4
\end{array}
\right],
\left[
\begin{array}{c}
  2\\
  5
\end{array}
\right],
\left[
\begin{array}{c}
  3\\
  6
\end{array}
\right].
\end{displaymath}

The investigation of these objects for greater dimensions shows the way
of the problem solution.


\subsection{Related problems}
The existence of the polynomial algorithm of graphs' isomorphism testing
allows to suggest the existence of a full polynomial invariant of a graph
and then the existence of a full polynomial invariant of any finite group,
because any finite group $G$ has a special for this group graph
representation as a partition of the set $G^2$ generated by the left
(right) action of group $G$ on the set $G^2$. Indeed, this graph can be
simplified to a $d+1$ color digraph (where one color is empty) with
geometrical structure of many-dimensional tore by dimension $d$ being
equal to a number of group's generators. The canonization of group's
generating system can allow to identify a group by full invariant of its
graph.

On the other hand, the existence of the polynomial algorithm of graphs'
isomorphism testing shows a possible existence of the full polynomial
invariant of a partition of $V^l$. It gives a different approach for group
identification as a permutation representation of minimal degree, acting
on $V$, and then as automorphism group of its system of orbits on $V^l$
for minimal possible $l$; where for symmetric group $l=1$, for graph
automorphism group $l=2$, for alternating group $l=n-1$.

It can be seen that described symmetries are equivalent to symmetry of the
set of left (right) cosets of cyclic subgroups of group $G$ generated by
the group's generators. This shows that considered approach is different
from the group symmetry investigation through the structure of the set of
its normal subgroups.

It should be said further if above considerations are correct then it
follows that an independent natural full invariant of an abstract group
is needed to exist. An investigation of invariants of described symmetrical
objects generated by these groups gives the way for obtaining this full
invariant.


\section*{Acknowledgements}
I would like to express many thanks  to Dr. M.Tabachnikov for the proposal
(in 1983) to solve this problem, to Dr.V.Grinberg for productive contacts,
to Dr. M.Klin for acquainting with contemporary achievements in this field
and to my sons Dr.~S.~Golubchik and R.~Golubchik for priceless computer
assistance.

\end{document}